\documentclass[11pt]{article}
\usepackage{amsfonts}
\usepackage{amsthm}

\begin{document}

\def\a{\alpha} 
 \def\b{\beta}
 \def\e{\epsilon}
 \def\d{\delta}
  \def\D{\Delta}
 \def\c{\chi}
 \def\k{\kappa}
 \def\g{\gamma}
 \def\t{\tau}
\def\ti{\tilde}
 \def\N{\mathbb N}
 \def\Q{\mathbb Q}
 \def\Z{\mathbb Z}
 \def\C{\mathbb C}
 \def\F{\mathbb F}
 \def\G{\Gamma}
 \def\go{\rightarrow}
 \def\do{\downarrow}
 \def\ra{\rangle}
 \def\la{\langle}
 \def\fix{{\rm fix}}
 \def\ind{{\rm ind}}
 \def\rfix{{\rm rfix}}
 \def\diam{{\rm diam}}
 \def\M{{\cal M}}
 
 \def\rank{{\rm rank}}
 \def\soc{{\rm soc}}
 \def\Cl{{\rm Cl}}
 \def\A{{\rm Alt}}
 
 \def\Sym{{\rm Sym}}
 \def\PSL{{\rm PSL}}
 \def\SL{{\rm SL}}
 \def\Ext{{\rm Ext}}
 \def\E{{\cal E}}
 \def\l{\lambda}
 \def\Lie{\rm Lie}
 \def\s{\sigma}
 \def\O{\Omega}
 \def\o{\omega}
 \def\ot{\otimes}
 \def\op{\oplus}
 \def\pf{\noindent {\bf Proof.$\;$ }}
 \def\Proof{{\it Proof. }$\;\;$}
 \def\no{\noindent}
 %{\quad\sqbox\vspace{1mm}\par}
\def\hal{\unskip\nobreak\hfil\penalty50\hskip10pt\hbox{}\nobreak
 \hfill\vrule height 5pt width 6pt depth 1pt\par\vskip 2mm}

 \renewcommand{\thefootnote}{}

\newtheorem{theorem}{Theorem}
 \newtheorem{thm}{Theorem}[section]
 \newtheorem{theor}{Theorem}[subsection]
 \newtheorem{prop}[thm]{Proposition}
 \newtheorem{lem}[theor]{Lemma}
 \newtheorem{lemma}[thm]{Lemma}
 \newtheorem{propn}[theor]{Proposition}
 \newtheorem{cor}[thm]{Corollary}
 \newtheorem{coroll}[theorem]{Corollary}
 \newtheorem{coro}[theor]{Corollary}
 \newtheorem{rem}[thm]{Remark}
 \newtheorem{cla}[thm]{Claim}

\parskip 1mm
  \topmargin=-20mm
    \textheight=237mm
    \textwidth=140mm

\title{Some word maps that are non-surjective on infinitely 
many finite simple groups}
\author{Sebastian Jambor \and Martin W. Liebeck \and E.A. O'Brien }
\date{}
\maketitle

\begin{abstract} 
We provide the first examples of words in the free group of rank~2 which 
are not proper powers and for which the 
corresponding word maps are non-surjective on an infinite family 
of finite non-abelian simple groups.
\end{abstract}

\section{Introduction}

The theory of word maps on finite non-abelian 
simple groups -- that is, maps of the form $(x_1,\ldots ,x_k)
\rightarrow w(x_1,\ldots ,x_k)$ for some word $w$ in the free group $F_k$ of
rank $k$ -- has attracted much recent attention. It was
shown in \cite[1.6]{LS} that for a given nontrivial word $w$, every
element of every sufficiently large finite simple group $G$ can be
expressed as a product of $C(w)$ values of $w$ in $G$, where $C(w)$
depends only on $w$; and this has been dramatically improved to $C(w)
= 2$ in \cite{LarSh, LarShT, Sh}. Improving $C(w)$ to 1 is not possible in
general, as is shown by power words $x_1^n$, which cannot be
surjective on any finite group of order non-coprime to $n$. 

Certain words are surjective on all groups -- namely, those in cosets of 
the form $x_1^{e_1}....x_k^{e_k}F_k'$ where the $e_i$ are integers with ${\rm gcd}(e_1,...,e_k)=1$
(see \cite[3.1.1]{S}). 
The word maps for a small number of other words have been shown to be surjective on all finite simple groups.
These include the commutator word $[x_1,x_2]$ (the Ore
conjecture \cite{ore}), the words $x_1^px_2^p$ (for a prime $p$) and
variants \cite{GM, LOST}. Other studies have restricted the simple
groups under consideration to families such as $\PSL_2(q)$ (see, 
for example, \cite{BGG}). Motivating some of this work is a conjecture of 
Shalev, stated in
\cite[Conjecture 8.3]{BGG}: 
if $w(x_1,x_2)$ is not a proper power of a non-trivial word, 
then the corresponding word map is surjective on $\PSL_2(q)$ for all
sufficiently large $q$.

Theorem \ref{word} gives a family of words which are counterexamples to Shalev's
conjecture. We believe these are the first non-power words
to be proved non-surjective on an infinite family of finite simple groups.

\begin{theorem} \label{word}
Let $k \geq 2$ be an integer such that $2k+1$ is prime, 
and let $w$ be the word $x_1^2[x_1^{-2}, x_2^{-1}]^k$.
Let $p \neq 2k + 1$ be a prime of inertia degree~$m > 1$ in 
$\Q(\zeta + \zeta^{-1})$, 
where $\zeta$ is a primitive $(2k+1)$-th root of unity, and 
$\left(\frac{2}{p}\right) = -1$.
Then the word map $(x,y) \to w(x,y)$ is non-surjective on $\PSL_2(q)$ for all 
$q = p^n$ where $n$ is a positive integer not divisible by $2$ or by~$m$.
\end{theorem}

%The following corollary is almost immediate from the theorem (see the end of the paper for the deduction).

\begin{coroll}\label{cor1}
Let $k \geq 2$ be an integer such that $2k+1$ is prime, 
and let $w$ be the word $x_1^2[x_1^{-2}, x_2^{-1}]^k$.
Let $p \ne 2k+1$ be an odd prime such that $p^2 \not\equiv 1 \bmod 16$ 
and $p^2 \not\equiv 1 \bmod (2k+1)$, and 
let $m$ be the smallest positive integer with 
$p^{2m} \equiv 1 \bmod (2k+1)$.
Then the word map $(x,y) \to w(x,y)$ is non-surjective on $\PSL_2(q)$ 
for all $q = p^n$ where 
$n$ is a positive integer not divisible by $2$ or by~$m$.
\end{coroll}

\noindent 
The corollary will be deduced from Theorem \ref{word} at the end of the paper. Taking $k=2$ we obtain the following. 

\begin{coroll} If $w = x_1^2[x_1^{-2}, x_2^{-1}]^2$, then the word map $(x,y) \to w(x,y)$ is non-surjective on $\PSL_2(p^{2r+1})$ for all non-negative integers $r$ and all odd primes~$p \neq 5$ such that $p^2 \not\equiv 1 \bmod 16$ and 
$p^2 \not\equiv 1 \bmod 5$.
\end{coroll}

\section{Proof of Theorem \ref{word}}
Let $K$ be a field and $G = \SL_2(K)$, and let $\c : G \rightarrow K$ be the trace map. 
A classical result of Fricke and Klein implies for every word $w \in F_2$, the free group of rank 2, 
there is a unique polynomial $\tau(w) \in \Z[s,t,u]$ such that 
for all $x,y \in G$, $\c(w(x,y))$ is equal to $\tau(w)$ evaluated at 
$s=\c(x)$, $t= \c(y)$, $u = \c(xy)$. We call $\tau(w)$ the {\it trace} polynomial of $w$. 
A proof of this fact, providing a 
constructive method of computing $\tau(w)$ for a given word~$w$, can
be found in \cite[2.2]{PF}. 
The method is based on the following identities for traces of $2\times 2$
matrices $A,B$ of determinant 1:
\[
\begin{array}{l}
{\rm Tr}(AB) = {\rm Tr}(BA) \\
{\rm Tr}(A^{-1}) = {\rm Tr}(A) \\
{\rm Tr}(A^2B) = {\rm Tr}(A){\rm Tr}(AB) - {\rm Tr}(B).
\end{array}
\]

\begin{lemma}
\label{L:prefactorization}
For $k \in \N$ and $w \in F_2$,
\[
    (-1)^k + \sum_{i = 1}^k (-1)^{k-i}\tau(w^i) = \prod_{i = 1}^k (\tau(w) + \zeta^i + \zeta^{-i}),
\]
where $\zeta$ is a primitive $(2k+1)$-th root of unity.
\end{lemma}

\pf 
We adapt the proof of \cite[Proposition~2.6]{PF}.
Assume first that $w = x_1$.
%For any $w, w_1, w_2 \in F_2$,
%\[
%    \tau(w(w_1, w_2)) = \tau(w)(\tau(w_1), \tau(w_2), \tau(w_1w_2)),
%\]
%so it is enough to prove the lemma for $w = x_1$.
Let $A := {\small \pmatrix {0 & 1 \cr -1 & s}}$.
By the uniqueness of the trace polynomial,
$\tau(w^i) = {\rm Tr}(A^i) = {\rm Tr} {\small \pmatrix {y^i & 0 \cr 0 & y^{-i}}} = y^i + y^{-i}$,
where $y + y^{-1} = s$. 
Hence
\[
\begin{array}{ll}
    \sum_{i = 1}^k (-1)^{k-i}\tau(w^i) + (-1)^k & = \sum_{i = 1}^k (-1)^{k-i}y^i + \sum_{i = 1}^k (-1)^{k-i}y^{-i} + (-1)^k \\
    & = y^{-k}\sum_{i = 0}^{2k} (-1)^i y^i \\
    & = y^{-k}\prod_{i = 1}^{2k} (y + \zeta^i) \\
    & = \prod_{i = 1}^k (y + \zeta^i)(1 + \zeta^{-i}y^{-1}) \\
    &    = \prod_{i = 1}^k (s + \zeta^i + \zeta^{-i}). 
\end{array}
\]
Note that for $v, v_1, v_2 \in F_2$,
\[
    \tau(v(v_1, v_2)) = \tau(v)(\tau(v_1), \tau(v_2), \tau(v_1v_2)),
\]
so the general case is derived from the special case $w = x_1$ by
polynomial evaluation at $s = \tau(w)$, i.e., setting $v = x_1^i$, $v_1 = w$, $v_2 = 1$.
\hal 

\begin{lemma}\label{l2}
Let $k \in \N$. The trace polynomial of 
$w = x_1^2[x_1^{-2},x_2^{-1}]^k$ factors over $\Z[\zeta + \zeta^{-1}]$ as
\[
    (s^2 - 2)\prod_{i = 1}^k (s^4 - s^3tu + s^2t^2 + s^2u^2 - 4s^2 + 2 + \zeta^i + \zeta^{-i}),
\]
where $\zeta$ is a primitive $(2k+1)$-th root of unity.
\end{lemma}

\pf
Let $c = [x_1^{-2},x_2^{-1}]$. We claim that 
\[
\tau (x_1^2c^k) = (\tau (x_1)^2 - 2) ( \sum_{i = 1}^k (-1)^{k-i}\tau(c^i) + (-1)^k).
\]
The result then follows by Lemma~\ref{L:prefactorization}, since $\tau(x_1) = s$ 
and $\tau(c) = s^4 - s^3tu + s^2t^2 + s^2u^2 - 4s^2 + 2$.

The proof is by induction on $k$.
The claim is easily verified for $k = 1, 2$.
For $k > 1$ it is equivalent to $\tau(x_1^2c^k) = (\tau(x_1)^2 - 2)\tau(c^k) - \tau(x_1^2c^{k-1})$.
Using the rule $\tau(x^2y) = \tau(x)\tau(xy) - \tau(y)$ for all $x,y \in F_2$ and 
the fact that $x_1^{-2}x_2^{-1} = x_2^{-1}x_1^{-2}c$, we deduce that
\[
\begin{array}{ll}
    \tau(x_1^2c^k) & = (\tau(x_1)^2 -1)\tau(c^k) - \tau(x_1)\tau(x_1x_2x_1^{-2}x_2^{-1}c^{k-1}) \\
    & = (\tau(x_1)^2 -1)\tau(c^k) - \tau(x_1)\tau(x_1^{-1}c^k) \\
    & = (\tau(x_1)^2 -1)\tau(c^k) - \tau(x_1^{-2}c^k) - \tau(c^k).
\end{array}
\]
Thus it suffices to prove that $\tau(x_1^{-2}c^k) = \tau(x_1^2c^{k-1})$.
Now $\tau(x_1^{-2}c^k) = \tau(c)\tau(c^{k-1}x_1^{-2}) - \tau(c^{k-2}x_1^{-2})$. By induction, for $k\ge 3$
this is equal to $\tau(c)\tau(x_1^2c^{k-2}) - \tau(x_1^2c^{k-3})$, which is 
equal to $\tau(x_1^2c^{k-1})$.
%This completes the proof of the lemma.
\hal

\vspace{2mm}
\no {\bf Proof of Theorem \ref{word}} \\
Let $q=p^n$ be as in the hypothesis of the theorem, let $K = \F_q$, and let $w$ be the word $x_1^2[x_1^{-2}, x_2^{-1}]^k$.
The ring of integers of $\Q(\zeta + \zeta^{-1})$ is $\Z[\zeta + \zeta^{-1}]$ (see \cite[Proposition~2.16]{washington}).
Since $2k+1$ is prime, $\Z[\zeta + \zeta^{-1}] = \Z[\zeta^i + \zeta^{-i}]$ for every $1 \leq i \leq k$.
Let $P \trianglelefteq \Z[\zeta^i + \zeta^{-i}]$ be a prime above $p$.
Then $\Z[\zeta^i + \zeta^{-i}]/P = \F_{p^m}$,
in particular $\zeta^i + \zeta^{-i}$ is a primitive element of $\F_{p^m}$ for every $1 \leq i \leq k$.

Suppose that some triple $(s,t,u) \in \F_q^3$ is a zero of the trace polynomial $\tau(w)$.
By Lemma~\ref{l2}, $\tau(w)$ factors as
\[
    (s^2 - 2)\prod_{i = 1}^k (s^4 - s^3tu + s^2t^2 + s^2u^2 - 4s^2 + 2 + \zeta^i + \zeta^{-i}),
\]
over $\F_{p^m}$, so $(s,t,u) \in \F_q^3 \subseteq \F_{q^m}^3$ must be a zero of one of the factors.
Since $s^2 - 2$ is irreducible over $\F_q$, $(s,t,u)$ must be a zero of 
$s^4 - s^3tu + s^2t^2 + s^2u^2 - 4s^2 + 2 + \zeta^i + \zeta^{-i}$ for some $i$.
This implies that $\zeta^i + \zeta^{-i} \in \F_q$, which is a 
contradiction.
Hence no element of $\SL_2(q)$ of the form $w(x,y)$ can have trace zero. 
%completing the proof of the theorem.
\hal

\vspace{2mm}
\no {\bf Proof of Corollary \ref{cor1} }  \\ 
Let $q = p^n$ be as in the hypothesis of the corollary.
The hypothesis $p^2 \not\equiv 1 \bmod 16$ is equivalent to $\left(\frac{2}{p}\right) = -1$.
By the cyclotomic reciprocity law (see for example \cite[Theorem~2.13]{washington}),
the inertia degree of $p$ in $\Q(\zeta)$ is $m$ or $2m$. In the former 
case, $m$ must be odd. 
Thus in both cases the inertia degree of $p$ in 
$\Q(\zeta + \zeta^{-1})$ is $m$, since $\Q(\zeta + \zeta^{-1})$
is a subfield of index~$2$ in $\Q(\zeta)$.
Now $p^2 \not\equiv 1 \bmod (2k+1)$ implies $m > 1$, 
and the conclusion follows from Theorem \ref{word}.
\hal

\vspace{0.5mm} 
\no {\bf Remark.} Our search for non-surjective words 
was assisted by \cite{CR}, which lists 
representatives of minimal length for certain automorphism
classes of words in $F_2$.

\vspace{2mm} 
\no {\bf Acknowledgements.} The first author was
supported by the DFG priority program SPP 1388. The second
author is grateful to the Royal Society for an International Exchange
grant, and to Nikolay Nikolov for very helpful discussions. All
authors were supported by the Marsden Fund of New Zealand
via grant UOA 1015.

\vspace{8mm}
\no Sebastian Jambor, Lehrstuhl B f\"ur Mathematik, RWTH Aachen, D-52056 Aachen, Germany.  Email: sebastian.jambor@rwth-aachen.de

\vspace{2mm}
\no Martin W.\ Liebeck, Department of Mathematics, Imperial College, London SW7 2AZ, UK.  Email: m.liebeck@imperial.ac.uk

\vspace{2mm}
\no E.A.\ O'Brien, Department of Mathematics, University of Auckland, New Zealand.  Email: e.obrien@auckland.ac.nz


\begin{thebibliography}{99}
\bibitem{BGG}
T. Bandman, S. Garion and F. Grunewald, On the surjectivity of Engel words 
on \PSL(2,q), to appear in {\it Groups, Geometry and Dynamics}.

\bibitem{CR}
B. Cooper and E. Rowland, 
Growing words in the free group on two generators,
to appear in {\it Illinois J. Math.}

\bibitem{GM} R.M. Guralnick and G. Malle, Products of conjugacy classes 
and fixed point spaces, {\it J. Amer. Math. Soc.} {\bf 25} (2012), 77--121.

\bibitem{LarSh} M. Larsen and A. Shalev, Word maps and Waring type problems. 
{\it J. Amer. Math. Soc.} {\bf 22} (2009), 437--466. 

\bibitem{LarShT} M. Larsen, A. Shalev and P.H. Tiep, The Waring problem for finite simple groups, {\it Annals of Math.} {\bf 174} (2011), 1885--1950.

\bibitem{LS} M.W. Liebeck and A. Shalev, Diameters of finite simple groups: 
sharp bounds and applications, {\it Annals of Math.} {\bf 154} (2001), 383--406. 

\bibitem{ore} M.W. Liebeck, E.A. O'Brien, A. Shalev, and P.H. Tiep, 
The Ore conjecture, {\it J. Eur. Math. Soc.} {\bf 12} (2010), 939--1008. 

\bibitem{LOST} 
M.W. Liebeck, E.A. O'Brien, A. Shalev, and P.H. Tiep, Products of squares 
in finite simple groups, 
{\it Proc. Amer. Math. Soc.} {\bf 140} (2012), 21--33. 

%\bibitem{neukirch} J. Neukirch, {\it Algebraic number theory}, Grundlehren der Mathematischen Wissenschaften, Vol. 322,
%   Springer-Verlag, Berlin, 1999.
%
\bibitem{PF} W. Plesken and A. Fabia\'nska, 
An $L_2$-quotient algorithm for finitely presented groups, 
{\it J. Algebra} {\bf 322} (2009), 914--935. 

\bibitem{S} D. Segal, {\it Words: notes on verbal width in groups}, London Math. Soc. Lecture Note Series 
{\bf 361}, Cambridge University Press, Cambridge, 2009.

\bibitem{Sh} A. Shalev, Word maps, conjugacy classes, and a noncommutative Waring-type theorem, 
{\it Annals of Math.} {\bf 170} (2009), 1383--1416.

\bibitem{washington} L.~C.\ Washington,
{\it Introduction to cyclotomic fields}, Graduate Texts in Mathematics, Vol.~83,
Springer-Verlag, New York, 1982.
\end{thebibliography}
\end{document}